\documentclass[11pt]{article} 
\usepackage{amsmath}
\usepackage{amsfonts}
\usepackage{graphicx}
\usepackage{color}
\usepackage[latin1]{inputenc}
\usepackage[T1]{fontenc}
\usepackage[french]{babel}
\makeatletter
\def\@sect#1#2#3#4#5#6[#7]#8{%
  \ifnum #2>\c@secnumdepth
    \let\@svsec\@empty
  \else
    \refstepcounter{#1}%
    \protected@edef\@svsec{\@seccntformat{#1}\relax}%
  \fi
  \@tempskipa #5\relax
  \ifdim \@tempskipa>\z@
    \begingroup
      #6{%
        \@hangfrom{\hskip #3\relax\@svsec}%
          \interlinepenalty \@M #8\@@par}%
    \endgroup
    \csname #1mark\endcsname{#7}%
    \addcontentsline{toc}{#1}{%
      \ifnum #2>\c@secnumdepth \else
        \protect\numberline{\csname the#1\endcsname.}%
      \fi
      #7}%
  \else
    \def\@svsechd{%
      #6{\hskip #3\relax
      \@svsec #8}%
      \csname #1mark\endcsname{#7}%
      \addcontentsline{toc}{#1}{%
        \ifnum #2>\c@secnumdepth \else
          \protect\numberline{\csname the#1\endcsname.}%
        \fi
        #7}}%
  \fi
  \@xsect{#5}}
\def\@seccntformat#1{\csname the#1\endcsname.\quad}
\makeatother
\newtheorem{theo}[equation]{Th\'eor\`eme}

\newtheorem{proposition}[equation]{Proposition}

\makeatletter

\@addtoreset{equation}{section}
\makeatother
\newcommand{\carrenoir}{\rule{0.5em}{0.5em}}
\makeatletter
\newenvironment{demo}[1][\@empty]{\textbf{D\'emonstration~%
\ifx\@empty#1:\else #1~:\fi~}}
{\hfill\carrenoir\nolinebreak\vspace{2mm}}
\makeatother

\newcommand{\oper}[2]{\newcommand{#1}{\mathop{\mathrm{#2}}\nolimits} }
\oper{\Vol}{Vol}
\newcommand{\R}{\mathbb R}
\newcommand{\de}{\mathrm{ d }}
\oper{\SO}{SO}

\title{Première valeur propre du laplacien, volume conforme et chirurgies}
\author{Pierre Jammes}
\date{}
\begin{document}
\maketitle
{\small 
\textsc{Résumé.---}
On définit dans cet article un nouvel invariant differentiel des 
variétés compactes par $V_{\mathcal M}(M)=\inf_g V_c(M,[g])$, où
$V_c(M,[g])$ désigne le volume conforme de la variété $M$ pour la 
classe conforme $[g]$, et on montre que cet invariant est uniformément
majoré. La principale motivation est qu'on en déduit une majoration
de l'invariant de Friedlander et Nadirashvili défini par
$\inf_g\sup_{\tilde g\in[g]}\lambda_1(M,\tilde g)\Vol(M,\tilde g)^{\frac 2n}$.

La démonstration consiste à étudier comment évolue l'invariant
$V_{\mathcal M}(M)$ quand on pratique des chirugies sur $M$.

Mots-clefs : première valeur propre du laplacien, volume conforme,
chirurgies.

\medskip
\textsc{Abstract.---}
We define a new differential invariant a compact manifold by
$V_{\mathcal M}(M)=\inf_g V_c(M,[g])$, where $V_c(M,[g])$ is the conformal
volume of $M$ for the conformal class $[g]$, and prove that it is uniformly
bounded above. The main motivation is that this bound provides a upper
bound of the Friedlander-Nadirashvili invariant defined by
$\inf_g\sup_{\tilde g\in[g]}\lambda_1(M,\tilde g)\Vol(M,\tilde g)^{\frac 2n}$.

The proof relies on the study of the behaviour of $V_{\mathcal M}(M)$
when one performs surgeries on $M$.


}

\section{Introduction}
 Si $(M^n,g)$ est une variété riemannienne compacte de dimension $n$, on 
notera $\lambda_1(M,g)$ la première valeur propre du laplacien sur $M$ pour
la métrique $g$.
 Dans \cite{fn99}, L.~Friedlander et N.~Nadirashvili ont défini un
nouvel invariant différentiable des variétés compactes en posant
\begin{equation}\label{intro:nu}
\nu(M)=\inf_g\sup_{\tilde g\in[g]}\lambda_1(M,\tilde g)
\Vol(M,\tilde g)^{\frac 2n},
\end{equation}
où $[g]=\{h^2g, h\in C^\infty, h>0\}$ désigne la classe conforme de 
la métrique riemannienne $g$.
L'invariant $\nu(M)$ est bien défini car on sait que $\sup_{\tilde g\in[g]}
\lambda_1(M,\tilde g)\Vol(M,\tilde g)^{\frac 2n}<+\infty$ (voir \cite{ly82},
\cite{esi86} et \cite{fn99}).

L'invariant $\nu(M)$ est très mal connu. En dimension~2, on peut déduire
de \cite{ly82} que $\nu(S^2)=8\pi$ et $\nu(\R P^2)=12\pi$ en utilisant
le fait que ces deux surfaces n'ont qu'une seule classe conforme, et
A.~Girouard a montré dans \cite{gi05} que $\nu(T^2)=\nu(K^2)=8\pi$, 
où $K^2$ désigne la bouteille de Klein. Mais en 
dimension plus grande, on sait seulement que $\nu(M^n)\geq\nu(S^n)=
n\Vol(S^n,g_{\mathrm{can}})^{\frac2n}$ (voir \cite{fn99} et \cite{ces03}).

 L'objet de cet article est de montrer qu'à dimension fixée, l'invariant
$\nu(M)$ est uniformément majoré :
\begin{theo}\label{intro:th1}
Il existe une constante $c(n)>0$ telle que pour toute variété compacte $M$
de dimension~$n$, on a $\nu(M)<c$.
\end{theo}

La démonstration du théorème~\ref{intro:th1} s'appuie essentiellement
sur des arguments géométriques et topologiques, ce qui contraste avec
les résultats d'A.~Girouard qui utilisent surtout des techniques
d'analyse. En effet,
le point de départ de la démonstration consiste à se ramener à un 
problème sur le volume conforme. Cet invariant conforme des variétés
compactes est défini par
\begin{equation}\label{intro:vc}
V_c(M,[g])=\inf_{\stackrel{N\geq 1}{\varphi:(M,[g])\hookrightarrow S^N}}
\sup_{\gamma\in G_N}\Vol(\gamma\circ\varphi(M)),
\end{equation}
l'application $\varphi$ parcourant l'ensemble des immersions
conformes de $(M,[g])$ dans $S^N$, et $G_N$ désignant le groupe de Möbius de 
dimension~$N$, c'est-à-dire le groupe des difféomorphismes conformes de $S^N$.
La propriété du volume conforme qui nous intéresse ici est qu'il permet
de majorer uniformément la première valeur propre du laplacien sur une
classe conforme:
\begin{theo}\label{intro:vcth}
Si $(M,g)$ est une variété riemannienne compacte, alors
$$\lambda_1(M,g)\Vol(M,g)^{\frac2n}\leq nV_\mathrm{c}(M,[g])^{\frac2n}.$$
\end{theo}
Ce théorème est démontré pour $n=2$ dans \cite{ly82}, et généralisé
en toute dimension par A.~El~Soufi et S.~Ilias dans \cite{esi86}.

On peut déduire immédiatement du théorème~\ref{intro:vcth} la majoration 
$\nu(M)\leq n\cdot\inf_gV_\mathrm{c}(M,[g])^{\frac2n}$. On est donc
ramené à majorer l'invariant différentiel $V_\mathcal{M}(M)=
\inf_gV_\mathrm{c}(M,[g])$, qu'on appellera \emph{volume de Möbius} dans
la suite de ce texte. Cette dénomination ce justifie par le fait qu'on
peut reformuler la définition de cet invariant par
\begin{equation}\label{intro:vm}
V_\mathcal{M}(M)=\inf_{\stackrel{N}{\varphi:M\hookrightarrow S^N}}
\sup_{\gamma\in G_N}\Vol(\gamma\circ\varphi(M)),
\end{equation}
la différence avec le volume conforme étant que $\varphi$ parcourt 
l'ensemble de toutes les immersions de $M$ dans $S^N$ sans restriction.
La référence aux métriques et au classes conformes sur $M$ a donc disparu,
seuls interviennent la métrique canonique de la sphère $S^N$ et l'action du
groupe de Möbius $G_N$. Notons aussi que la minoration classique 
$V_\mathrm{c}(M^n,[g])\geq\Vol(S^n,g_\textrm{can})$ implique que
$V_\mathcal{M}(M^n)\geq\Vol(S^n,g_\textrm{can})$.

On va donc montrer le
\begin{theo}\label{intro:th2}
Il existe une constante $c(n)>0$ telle que pour toute variété compacte $M$ 
de dimension~$n$, on a $V_\mathcal{M}(M)<c$.
\end{theo}
L'idée de la démonstration consiste à étudier comment varie le volume de 
Möbius quand on pratique des chirurgies sur la variété. Nous distinguerons
deux cas: dans la section~\ref{2} nous traiterons le cas de la dimension~2
en expliquant les aspects géométriques de la démonstration, et dans
la section~\ref{n} nous verrons le cas des dimensions plus grandes, pour
lequel les aspects topologiques demandent plus de travail.

 Notons qu'en dimension~2, le théorème~\ref{intro:th2} peut se déduire
des travaux sur les immersions de surfaces dans $\R^3$. On définit la 
fonctionnelle de Willmore d'une surface $\Sigma$ immergée dans $\R^3$
par $\mathcal W(\Sigma)=\int_\Sigma H^2\de v$, où $H$ désigne la courbure
moyenne de $\Sigma$. D'une part, P.~Li et S.~T.~Yau obtiennent dans 
\cite{ly82} l'inégalité $V_\mathrm{c}(\Sigma)\leq\mathcal W(\Sigma)$ et 
d'autre part R.~Kusner a montré dans \cite{ku89} que toute surface 
admet une immersion telle que $\mathcal W(\Sigma)<16\pi$. Il en découle
immédiatement que dans le théorème~\ref{intro:th2}, on peut choisir 
$c(2)<16\pi$.

 On peut se demander si l'égalité $V_\mathcal{M}(M^n)=\Vol(S^n,g_\textrm{can})$
caractérise la sphère canonique; rappelons que la même question a été posée 
au sujet du volume conforme par P.~Li et S.~T.~Yau dans \cite{ly82} et
qu'elle est encore ouverte actuellement.

Je remercie Bernd Ammann pour plusieurs discussions instructives autour
de la théorie du cobordisme.

\section{Dimension 2}\label{2}

Nous allons d'abord reformuler la définition du volume de Möbius.
Si $\gamma$ est un élément quelconque de $G_N$ et 
$r$ une isométrie de la sphère $S^N$, alors $\Vol(r\circ\gamma\circ\varphi
(M))=\Vol(\gamma\circ\varphi(M))$ pour toute immersion $\varphi$ de $M$
dans $S^N$. On peut donc écrire 
\begin{equation}\label{2:vm}
V_{\mathcal M}(M)=\inf_{\varphi:M\hookrightarrow S^N}
\sup_{\gamma\in G'_N}\Vol(\gamma\circ\varphi(M)),
\end{equation}
où $G'_N$ est une sous-ensemble de $G_N$ (pas nécessairement un 
sous-groupe) tel que tout élément de $G_N$ s'écrive sous la forme
$r\circ\gamma$ avec $\gamma\in G'_N$ et $r\in\SO(N)$.

Il existe plusieurs choix possibles pour $G'_N$, mais l'un d'entre
eux est particulièrement adapté à la démonstration qui suit. En projetant
stéréographiquement la sphère $S^N$ sur $\R^N\cup\{\infty\}$, on se
ramène à considérer des immersion $\varphi:M\hookrightarrow 
\R^N\cup\{\infty\}$ en calculant les volumes pour la métrique 
$g_{S^N}=\frac4{(1+\|x\|^2)^2}g_{\mathrm{eucl}}$ où $g_{\mathrm{eucl}}$
désigne la métrique euclidienne canonique. On peut alors choisir
pour $G'_N$ l'ensemble des homothéties et des translations de
$\R^N$ (voir le théorème~3.5.1 de~\cite{be83}).

 Supposons maintenant que $n=2$. On va montrer que pour une surface
compacte $M$ donnée, on peut contrôler le volume de Möbius de la
variété $\widetilde M$ obtenue par adjonction d'une anse à $M$ en
fonction du volume de Möbius de $M$. On pourra alors conclure par une 
récurrence sur le genre de la surface.

On commence par fixer un réel $\varepsilon>0$ et choisir une immersion 
$\varphi$ de $M$ dans un espace $\R^N$ telle que
\begin{equation}\label{2:extremal}
\sup_{\gamma\in G'_N}\Vol(\gamma\circ\varphi(M))-\varepsilon\leq
\Vol(\varphi(M))\leq V_\mathcal{M}(M)+\varepsilon,
\end{equation}
c'est-à-dire que $\varphi$ est proche de réaliser à la fois la borne
supérieure et la borne inférieure dans la définition~(\ref{2:vm}). On
va construire une immersion $\widetilde\varphi$ de $\widetilde M$ à partir de 
$\varphi$ et chercher à contrôler $\Vol(\gamma\circ\widetilde
\varphi(\widetilde M))$ pour tout $\gamma\in G'_N$.

\begin{figure}[h]
\begin{center}
\begin{picture}(0,0)%
\includegraphics{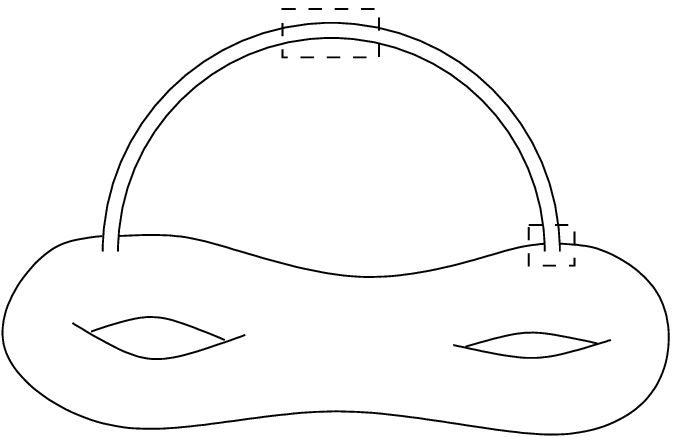}%
\end{picture}%
\setlength{\unitlength}{4144sp}%
\begingroup\makeatletter\ifx\SetFigFontNFSS\undefined%
\gdef\SetFigFontNFSS#1#2#3#4#5{%
  \reset@font\fontsize{#1}{#2pt}%
  \fontfamily{#3}\fontseries{#4}\fontshape{#5}%
  \selectfont}%
\fi\endgroup%
\begin{picture}(3070,1969)(923,-2835)
\put(3562,-1876){\makebox(0,0)[lb]{\smash{{\SetFigFontNFSS{9}{10.8}{\rmdefault}{\mddefault}{\updefault}{\color[rgb]{0,0,0}(b)}%
}}}}
\put(2377,-1249){\makebox(0,0)[lb]{\smash{{\SetFigFontNFSS{9}{10.8}{\rmdefault}{\mddefault}{\updefault}{\color[rgb]{0,0,0}(a)}%
}}}}
\end{picture}%
\end{center}
\caption{\label{2:surface}}
\end{figure}

La variété $\widetilde M$ est obtenue à partir de $M$ par adjonction d'une 
anse. On obtient donc une immersion $\widetilde\varphi$ de $\widetilde M$ 
en ajoutant une anse fine à $\varphi(M)$ (voir figure~\ref{2:surface}). 
Dans la suite, cette anse sera un tube centré sur un arc de courbe fixée 
et de rayon $\delta>0$ variable.
On peut en outre choisir cette courbe de sorte que l'anse ne coupe jamais
$\varphi(M)$, quitte à choisir $N$ suffisamment grand.

On doit chercher à contrôler le volume de $\gamma\circ\widetilde
\varphi(\widetilde M)$ pour tout $\gamma\in G'_N$. Pour cela, un élément 
crucial est le fait que si l'anse ou $\varphi(M)$ est projeté près du 
point $\infty$ sous l'action de $\gamma$, son volume devient négligeable 
pour la métrique $g_{S^N}$, même si le rapport d'homothétie de 
$\gamma$ est grand: en effet, une homothétie de rapport $R$ va multiplier
l'aire euclidienne par $R^2$, mais dans la formule 
$g_{S^N}=\frac4{(1+\|x\|^2)^2}g_{\mathrm{eucl}}$ la norme de $x$ est de 
l'ordre de $R$ pour les parties de $\widetilde\varphi(\widetilde M)$
projeté près de l'infini et donc leur aire devient petite quand $R$
devient grand.

Pour tout $\gamma\in G'_N$, on notera $R_\gamma$ son rapport d'homothétie.
Si on se donne une constante $R_0>0$, on peut trouver $\delta_0>0$ tel 
que si $\delta<\delta_0$ alors pour tout $\gamma$ tel que $R_\gamma\leq R_0$
le volume de l'anse transportée par $\gamma$ reste négligeable, c'est-à-dire 
que $\Vol(\gamma\circ\widetilde\varphi(\widetilde M))\leq
\Vol(\gamma\circ\varphi(M))
+\varepsilon$ et donc $\Vol(\gamma\circ\widetilde\varphi(\widetilde M))\leq 
V_\mathcal{M}(M)+3\varepsilon$ d'après (\ref{2:vm}). 
Si $R_0$ est suffisamment grand, l'argument précédent s'applique aussi 
pour $R_\gamma>R_0$ quand l'homothétie $\gamma$ envoie l'anse près de
l'infini, auquel cas son volume est aussi négligeable. Il reste donc
à traiter la situation où $R_\gamma>R_0$ et qu'une partie de l'anse reste 
proche de l'origine, qui se scinde en deux cas: $\varphi(M)$ est
envoyé près de l'infini par $\gamma$ (c'est-à-dire qu'on fait un zoom
sur (a) dans la figure~\ref{2:surface}), ou bien la la partie de
l'anse proche de l'origine est aussi proche du point d'attache à $\varphi(M)$ 
(zoom sur (b) dans la figure~\ref{2:surface}). 

\begin{figure}[h]
\begin{center}
\begin{picture}(0,0)%
\includegraphics{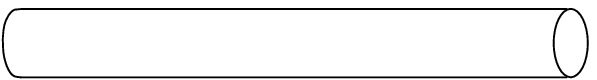}%
\end{picture}%
\setlength{\unitlength}{4144sp}%
\begingroup\makeatletter\ifx\SetFigFontNFSS\undefined%
\gdef\SetFigFontNFSS#1#2#3#4#5{%
  \reset@font\fontsize{#1}{#2pt}%
  \fontfamily{#3}\fontseries{#4}\fontshape{#5}%
  \selectfont}%
\fi\endgroup%
\begin{picture}(2695,337)(1216,-147)
\end{picture}%
\end{center}
\caption{\label{2:cylindre}}
\end{figure}

Dans le premier cas, on va comparer $\Vol(\gamma\circ\widetilde
\varphi(\widetilde M))$ au volume d'un cylindre infini 
(figure~\ref{2:cylindre}). Notons $\mathcal C$ le cylindre d'équation 
$x_1^2+x_2^2=1$ et $x_4=\ldots=x_N=0$ dans $\R^N$. Si $R_0$ est
suffisamment grand par rapport à la courbure de l'anse, le volume de 
$\gamma\circ\widetilde\varphi(\widetilde M)$ peut être approché par celui
d'un cylindre homothétique à $\mathcal C$, c'est-à-dire que
\begin{equation}
\Vol(\gamma\circ\widetilde\varphi(\widetilde M))\leq V_{\mathcal M}(\mathcal C)+
\varepsilon,
\end{equation}
où $V_{\mathcal M}(\mathcal C)=\sup_{\gamma'\in G'_N}
\Vol(\gamma'(\mathcal C))$. 

\begin{figure}[h]
\begin{center}
\begin{minipage}[b]{.46\linewidth}
\begin{picture}(0,0)%
\includegraphics{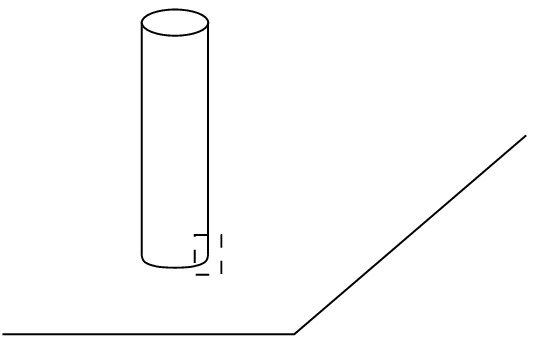}%
\end{picture}%
\setlength{\unitlength}{4144sp}%
\begingroup\makeatletter\ifx\SetFigFontNFSS\undefined%
\gdef\SetFigFontNFSS#1#2#3#4#5{%
  \reset@font\fontsize{#1}{#2pt}%
  \fontfamily{#3}\fontseries{#4}\fontshape{#5}%
  \selectfont}%
\fi\endgroup%
\begin{picture}(2417,1504)(1114,-1153)
\put(2156,-869){\makebox(0,0)[lb]{\smash{{\SetFigFontNFSS{8}{9.6}{\rmdefault}{\mddefault}{\updefault}{\color[rgb]{0,0,0}(c)}%
}}}}
\end{picture}%
\caption{\label{2:plan}}
\end{minipage}\hfill
\begin{minipage}[b]{.46\linewidth}
\begin{picture}(0,0)%
\includegraphics{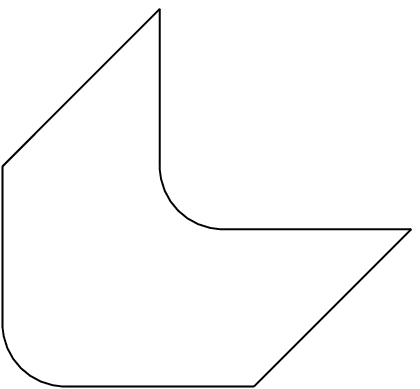}%
\end{picture}%
\setlength{\unitlength}{4144sp}%
\begingroup\makeatletter\ifx\SetFigFontNFSS\undefined%
\gdef\SetFigFontNFSS#1#2#3#4#5{%
  \reset@font\fontsize{#1}{#2pt}%
  \fontfamily{#3}\fontseries{#4}\fontshape{#5}%
  \selectfont}%
\fi\endgroup%
\begin{picture}(1893,1750)(1564,-1349)
\end{picture}%
\caption{\label{2:angle}}
\end{minipage}
\end{center}
\end{figure}

Dans le second cas, on peut choisir $R_0$ suffisamment grand par
rapport à la courbure de l'anse et de $\varphi(M)$, de sorte qu'on puisse 
approcher le volume de $\gamma\circ\widetilde\varphi(\widetilde M)$ par celui d'un plan
privé d'un disque auquel on recolle un demi-cylindre (figure~\ref{2:plan}).
On sait que le volume de ce plan pour la métrique sphérique est majoré
par le volume de la sphère $S^2$ (voir \cite{ly82} ou \cite{esi86}).
On a donc dans ce cas la majoration 
\begin{equation}
\Vol(\gamma\circ\widetilde\varphi(\widetilde M))\leq 4\pi+
V_{\mathcal M}(\mathcal C)+\varepsilon.
\end{equation}

Dans le cas précédent et dans la figure~\ref{2:plan}, on a considéré que
l'immersion $\widetilde\varphi$ n'est pas lisse au niveau de la jonction
de l'anse. Étant donné un réel 
$R_1\gg R_0$, on peut lisser $\widetilde\varphi$ de sorte que les 
approximations faites précédemment restent valables tant que 
$R_\gamma\leq R_1$. On choisit $R_1$ suffisamment grand pour que
si $R_\gamma>R_1$, la courbure du cylindre devienne négligeable.
La zone (c) dans la figure~\ref{2:plan} peut alors être représenté comme
sur la figure~\ref{2:angle} par deux demi-plans reliés par un quart de 
cylindre. le volume de $\gamma\circ\widetilde\varphi(\widetilde M)$ peut alors
être grossièrement majoré celui de deux plans et d'un cylindre, c'est-à-dire :
\begin{equation}
\Vol(\gamma\circ\widetilde\varphi(\widetilde M))\leq 8\pi+
V_{\mathcal M}(\mathcal C)+\varepsilon.
\end{equation}

Finalement, quel que soit $\gamma$, on a
\begin{equation}
\Vol(\gamma\circ\widetilde\varphi(\widetilde M))\leq 
\sup\left\{V_\mathcal{M}(M),8\pi+V_{\mathcal M}(\mathcal C)\right\}
+3\varepsilon,
\end{equation}
et ce pour tout $\varepsilon>0$. Si on pose $c=8\pi+
V_{\mathcal M}(\mathcal C)$, on a donc 
\begin{equation}
\sup_{\gamma\in G'_N}\Vol(\gamma\circ\widetilde\varphi(\widetilde M))\leq 
\sup\left\{V_\mathcal{M}(M),c\right\}
\end{equation}
et par conséquent $V_\mathcal{M}(\widetilde M)\leq\sup\left\{
V_\mathcal{M}(M),c\right\}$.

 Comme toutes les surfaces compactes peuvent s'obtenir par adjonction
d'anses à partir de la sphère et du plan projectif dont les volumes
de Möbius sont inférieurs à $c$, on en déduit
par une récurrence sur le genre de la surface que 
\begin{equation} V_\mathcal{M}(M)\leq c\end{equation}
pour toute surface compacte $M$.

\section{Dimension supérieure ou égale à 3}\label{n}
La démonstration du théorème~\ref{intro:th2} dans le cas général se déroule 
de la même manière qu'en dimension~2. L'idée est d'étudier comment évolue
de volume de Möbius quand on procède sur la variété à des chirurgies,
ce qui généralise l'adjonction d'anse. On se référera aux livres \cite{ra02}
(chapitres~1, 2 et 6) et \cite{ko07} (chapitre~VI et VII) pour les
aspects topologiques de la démonstration.

Plus précisément, si $S^k\hookrightarrow M^n$ est une sphère plongée
dans $M$ dont le fibré normal est trivial, on peut trouver un voisinage 
tubulaire de $S^k$ difféomorphe à $S^k\times B^{n-k}$. Une chirurgie 
de le long de $S^k$ consiste alors à construire une nouvelle variété
$\widetilde M$ en privant $M$ de ce voisinage tubulaire 
$S^k\times B^{n-k}$ et en recollant le long du bord ainsi créé la variété
$B^{k+1}\times S^{n-k-1}$, ce procédé pouvant s'effectuer de manière lisse.
On appelle degré de la chirurgie la dimension $k$ de la sphère $S^k$. Par 
exemple, l'adjonction d'anse vue dans la partie précédente est une chirurgie
de degré $0$ (voir le premier chapitre de \cite{ra02}).

Sur le modèle du paragraphe précédent, on obtient le résultat suivant :
\begin{proposition}\label{n:prop1}
Soit $n\geq3$ et $0\leq k\leq n-2$. Il existe une constante $c(n,k)>0$ telle
que si $\widetilde M$ est une variété obtenue par chirurgie de degré
$k$ sur une variété compacte $M$ de dimension $n$, alors
$V_\mathcal{M}(\widetilde M)\leq\sup\left\{
V_\mathcal{M}(M),c\right\}$.
\end{proposition}
Nous ne détaillerons pas la démonstration car c'est la même qu'en dimension~2:
après avoir choisi une sphère $S^k\hookrightarrow M$ à fibré normal 
trivial et une immersion 
$\varphi(M)$ vérifiant (\ref{2:extremal}), on obtient une immersion 
$\tilde\varphi(\widetilde M)$ en enlevant à $\varphi(M)$ un voisinage 
tubulaire de $\varphi(S^k)$ de rayon $\delta$ et en recollant le long
du bord une anse $B^{k+1}\times S^{n-k-1}(\delta)$. Le rôle du cylindre
$\mathcal C$ est joué par le produit $S^{n-k-1}\times\R^{k+1}$.
On voit que pour que le volume de l'anse tende vers zéro quand $\delta\to0$,
il est indispensable que la sphère $S^{n-k-1}$ soit de dimension au moins~1.
C'est la raison pour laquelle on a exclu le cas $k=n-1$ dans la 
proposition~\ref{n:prop1}.

 Des arguments classiques de la théorie du cobordisme permettent maintenant
de déterminer quelles variétés on peut obtenir à partir d'une variété donnée
à l'aide de chirurgies de degré au plus $n-2$. Rappelons que deux variétés 
compactes (pas nécessairement connexes) $M_1$ et $M_2$ de dimension $n$ 
sont cobordantes s'il existe une variété à bord $W$ de dimension $n+1$ 
dont le bord est l'union disjointe de $M_1\sqcup M_2$ et que la relation 
ainsi définie est une relation d'équivalence sur l'ensemble des variétés 
compactes. Plus précisément, on appelle \emph{cobordisme non orienté}
cette relation. Nous auront besoin en outre d'une autre notion de cobordisme, 
appelée \emph{cobordisme orienté} définie uniquement sur les variétés
orientées : une orientation sur la variété à bord $W$ induit une orientation 
sur son bord, deux variétés orientées $M_1$ et $M_2$ de dimension~$n$
sont alors cobordantes s'il existe une variété orientée $W$ dont le bord
est la réunion de $M_1$ et $-M_2$, $-M$ désignant la variété $M$
munie de son orientation opposée (on considère sur $\partial W$ l'orientation
induite par celle de $W$ et la donnée d'une normale sortante). 
Rappelons aussi que deux variétés orientées
peuvent être cobordantes au sens du cobordisme non orienté sans l'être
au sens du cobordisme orienté. Dans la suite de ce texte, les cobordismes
entre variétés orientables seront toujours considérés au sens orienté.

 La démonstration du théorème~\ref{intro:th2} va se poursuivre en traitant
de manière analogue, mais distincte, les cas des variétés orientables
et non orientables.
\begin{proposition}\label{n:prop2}
Soit $M_1$ et $M_2$ deux variétés connexes compactes de même dimension 
$n\geq3$. Si $M_1$ et $M_2$ sont cobordantes et toutes deux orientables 
(resp. non orientables), alors $M_2$ peut être obtenue à partir de $M_1$ 
par une suite finie de chirurgies de degré au plus $n-2$.
\end{proposition}
\begin{demo}
Rappelons d'abord quelques résultats classiques de la théorie de Morse de 
de la théorie du cobordisme (sur ce sujet, on peut consulter le chapitre~2
de \cite{ra02} ou le chapitre~VII de \cite{ko07}). 
Soit $W^{n+1}$ une variété dont le bord est $M_1\sqcup M_2$ (ou
$M_1\sqcup-M_2$ dans le cas orienté). On
peut trouver une fonction de Morse $f:W\to[0,1]$ telle que $f^{-1}(0)=M_1$,
$f^{-1}(1)=M_2$, dont les points critiques $x_1,\ldots, x_k$ vérifient 
$0<f(x_1)<\ldots<f(x_k)<1$ et telle que la suite $\lambda_i$ des indices des 
points critiques $x_i$ soit croissante (voir \cite{ko07}, ch.~VII, sections~1
et 2). On sait alors que si $y\in[0,1]$, les courbes de niveau $f^{-1}(y)$ 
sont difféomorphes entre elles quand $y$ varie entre deux valeurs critiques 
successives, et que si l'intervalle $[y,y']$ ne contient qu'une seule valeur 
critique $x_i$, on peut passer de $f^{-1}(y)$ à $f^{-1}(y')$ par une chirurgie 
de degré $\lambda_i-1$. On peut donc passer de $M_1$ à $M_2$ par une suite
finie de chirurgies de degrés croissants, et on veut montrer qu'on peut
se passer des chirurgies de degré $n-1$.

Si $W$ est orientée, la donnée d'un champ de vecteur transverse aux
courbes de niveau (par exemple le gradient de $f$) induit une orientation
sur ces courbes de niveau. Si $M_1$ et $M_2$ sont orientées, on peut donc 
se ramener au cas où on passe de $M_1$ à $M_2$ uniquement par des chirurgies
de degré $n-1$. On a la même conclusion si $M_1$ et $M_2$ sont non orientables
car une variété non orientable reste non orientable après des chirurgies de
degré 0 à $n-2$. On va montrer que dans les deux cas, ces chirurgies de 
degré~$n-1$ peuvent être remplacées par des chirurgies de degré~1.

\begin{figure}[h]
\begin{center}
\begin{picture}(0,0)%
\includegraphics{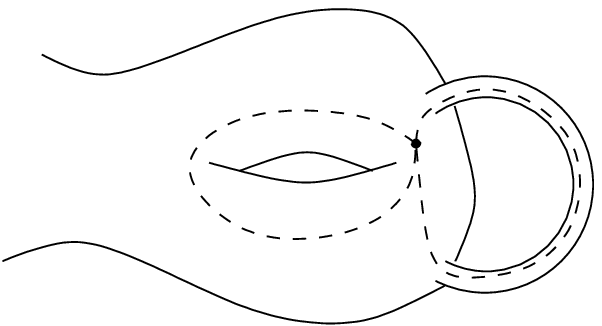}%
\end{picture}%
\setlength{\unitlength}{4144sp}%
\begingroup\makeatletter\ifx\SetFigFontNFSS\undefined%
\gdef\SetFigFontNFSS#1#2#3#4#5{%
  \reset@font\fontsize{#1}{#2pt}%
  \fontfamily{#3}\fontseries{#4}\fontshape{#5}%
  \selectfont}%
\fi\endgroup%
\begin{picture}(2721,1461)(1339,-1843)
\put(3306,-1331){\makebox(0,0)[lb]{\smash{{\SetFigFontNFSS{12}{14.4}{\rmdefault}{\mddefault}{\updefault}{\color[rgb]{0,0,0}$\gamma_1$}%
}}}}
\put(2366,-806){\makebox(0,0)[lb]{\smash{{\SetFigFontNFSS{12}{14.4}{\rmdefault}{\mddefault}{\updefault}{\color[rgb]{0,0,0}$\gamma_2$}%
}}}}
\end{picture}%
\end{center}
\caption{\label{n:chirurgie}}
\end{figure}
 Remarquons d'abord que l'opération inverse d'une chirurgie de degré $k$ 
est une chirurgie de degré $n-k-1$. En partant de $M_2$, on peut donc
obtenir $M_1$ en pratiquant des chirurgies de degré~0, c'est-à-dire
en ajoutant des anses. Ces chirurgies 
peuvent \emph{a priori} être de deux types, orientable (comme celle
qui permet de passer de $S^2$ à $T^2$) ou non orientable (comme pour
passer de $S^2$ à $K^2$, ce cas étant exclu si $M_1$ et $M_2$ sont
orientables). Dans le premier cas, on peut détruire l'anse
à l'aide du lemme d'élimination de Smale (voir \cite{ko07}, ch.~VI,
théorème~7.4)
en pratiquant une chirurgie le long d'un cercle $\gamma_1$ traversant 
l'anse comme sur la figure~\ref{n:chirurgie}. Cela suffit pour conclure
dans le cas où $M_1$ et $M_2$ sont orientables. Dans le cas où $M_1$ et $M_2$
ne sont pas orientables, le lacet $\gamma_1$ n'est pas nécessairement 
orientable (c'est-à-dire que son fibré normal n'est pas forcément trivial 
et par conséquent son voisinage tubulaire n'est pas un produit 
$S^1\times B^{n-1}$) et ne permet donc pas forcément de chirurgie. Mais 
comme $M_2$ est non orientable, on peut trouver 
un second lacet non orientable $\gamma_2$ dans $M_1$ qui est indépendant de 
l'anse, l'un des lacets $\gamma_1$ et $\gamma_1+\gamma_2$ est alors 
orientable et permet une chirurgie qui supprime l'anse.
\end{demo}

 On déduit aisément des proposition \ref{n:prop1} et \ref{n:prop2} que
le volume de Möbius est uniformément majoré sur une classe de cobordisme
orienté et sur l'ensemble des variétés non orientables d'une
classe de cobordisme non orienté.

Il reste à montrer qu'à dimension fixée, on peut trouver dans les deux cas
un majorant indépendant de la classe de cobordisme. Le cas non
orientable est le plus simple : R.~Thom a entièrement déterminé dans 
\cite{th54} l'ensemble des classes de cobordisme non orienté
(théorème~IV.12). En particulier, il n'y en a qu'un nombre fini à dimension
fixée.

Dans le cas orientable, on utilisera le fait que l'ensemble $\Omega_n$ des 
classes de cobordisme de dimension $n$ peut être muni d'une structure de 
groupe abélien en posant $M_1+M_2=M_1\#M_2$, l'élément neutre étant la 
classe de la sphère. R.~Thom a aussi montré dans \cite{th54} que ce
groupe est de type fini (théorème~IV.15), il suffit donc de majorer le
volume de Möbius d'une somme de variétés.

\begin{proposition}\label{n:prop3}
Il existe une constante $c(n)>0$ telle que si $M_1$ et $M_2$ sont deux 
variétés compactes de même dimension $n\geq3$ alors
$V_{\mathcal M}(M_1\#M_2)\leq\sup(V_{\mathcal M}(M_1),
V_{\mathcal M}(M_2),c)$.
\end{proposition}
\begin{demo}
On se donne deux immersions $\varphi_i$, $i=1,2$ de $M_i$ dans une même 
espace $\R^N$ muni de la métrique sphérique, pour un entier $N$ 
suffisamment grand, telles que
\begin{equation}\label{n:extremal}
\sup_{\gamma\in G'_N}\Vol(\gamma\circ\varphi_i(M_i))-\varepsilon\leq
\Vol(\varphi_i(M_i))\leq V_\mathcal{M}(M_i)+\varepsilon.
\end{equation}
 On se donne par ailleurs une translation $t_{\vec{v}}$ de vecteur $\vec v$,
et on définit une immersion $\varphi:M_1\sqcup M_2\to \R^N$ en posant
$\varphi_{|M_1}=\varphi_1$ et $\varphi_{|M_2}=t_{\vec v}\circ\varphi_2$.

Si $\vec v$ est suffisamment grand, le volume de $\varphi_{|M_2}$ est
négligeable et on a $\Vol(\varphi(M))\leq\Vol(\varphi_1(M_1))+\varepsilon$.
Si on se donne un élément $\gamma\in G'_N$ on a 
$\Vol(\gamma\circ\varphi(M_1))\leq\Vol(\varphi_1(M_1))+\varepsilon$,
et tant que $\gamma\circ\varphi(M_2)$ reste près de l'infini, son volume est
négligeable.

 Si $\gamma\circ\varphi(M_2)$ se rapproche de l'origine, deux situations
peuvent se produire. La première est que $\gamma\circ\varphi(M_1)$ est
envoyé près de l'infini, et on a alors 
\begin{equation}
\Vol(\gamma\circ\varphi(M_1\sqcup M_2))\leq\Vol(\gamma\circ\varphi(M_2))+
\varepsilon\leq V_\mathcal{M}(M_2)+3\varepsilon.
\end{equation}
La seconde possibilité est que $\gamma\circ\varphi(M_1)$ et 
$\gamma\circ\varphi(M_2)$ soient tous les deux proches de l'origine, 
mais alors le rapport d'homothétie de $\gamma$ est petit et le volume des
deux immersions est négligeable. 

 Dans tous les cas, on a bien la majoration 
$\Vol(\gamma\circ\varphi(M_1\sqcup M_2))\leq
\sup(V_{\mathcal M}(M_1),V_{\mathcal M}(M_2))+3\varepsilon$. Comme on peut
passer de $M_1\sqcup M_2$ à $M_1\#M_2$ par une chirurgie de degré~0, le
volume de Möbius de $M_1\#M_2$ est bien majoré.
\end{demo}

 À partir d'une majoration du volume de Möbius sur des représentants des
générateurs du groupe de cobordisme orienté, on obtient donc cette majoration
sur des représentants de toutes les classes de cobordisme.

\noindent Pierre \textsc{Jammes}\\
Université d'Avignon et des pays de Vaucluse\\
Laboratoire d'analyse nonlinéaire et géométrie (EA 2151)\\
F-84018 Avignon\\
\texttt{Pierre.Jammes@univ-avignon.fr}
\end{document}